\documentclass[times]{iapress}

\usepackage{moreverb}

\usepackage[colorlinks,bookmarksopen,bookmarksnumbered,citecolor=red,urlcolor=red]{hyperref} 

\def\volumeyear{2016}
\def\volumenumber{4}
\def\volumemonth{March}
\usepackage{amsmath,amssymb}
\usepackage{graphicx}

\newtheorem{example}{Example}
\newtheorem{lemma}{Lemma}
\newtheorem{theorem}{Theorem}
\newtheorem{corollary}{Corollary}
\newtheorem{definition}{Definition}

\newcommand{\lp}[1]{\left( \begin{array}{#1} }
\newcommand{\rp}{\end{array} \right)}

\newcommand{\me}{\mathbf}
\newcommand{\mr}{\mathbb}
\newcommand{\mt}{\mathsf}
\newcommand{\md}{\mathcal}
\newcommand{\ld}{\left}
\newcommand{\rd}{\right}
\newcommand{\be}{\begin{equation}}
\newcommand{\ee}{\end{equation}}
\newcommand{\ip}{\int_{-\pi}^{\pi}}
\newcommand{\eps}{\varepsilon}

\begin{document}

\runningheads{Filtering Problem for Functionals of Stationary Sequences}{M. Luz and M. Moklyachuk}

\title{Filtering Problem for Functionals of Stationary Sequences }

\author{Maksym Luz \affil{1},
    Mikhail Moklyachuk \affil{1}$^,$\corrauth
 }

\address{
\affilnum{1}Department of Probability Theory, Statistics and Actuarial
Mathematics, Taras Shevchenko National University of Kyiv, Kyiv, Ukraine
}

\corraddr{Mikhail Moklyachuk (Email: Moklyachuk@gmail.com). Department of Probability Theory, Statistics and Actuarial
Mathematics, Taras Shevchenko National University of Kyiv, Volodymyrska 64 Str., Kyiv 01601, Ukraine.}

\begin{abstract}
The problem of the mean-square optimal linear estimation of functionals
which depend on the unknown values of a stationary stochastic sequence from observations of the sequence with noise is considered. In the case of spectral certainty, where the spectral densities of the sequences are exactly known, we propose formulas for calculating the spectral characteristic and value of the mean-square error of the estimate, which are determined using the Fourier coefficients of some functions from the spectral densities. The minimax-robust method of estimation is applied in the case of spectral uncertainty, where the spectral densities are not exactly known, but a class of admissible spectral densities is given.
Formulas for determining the least favorable spectral densities and the
minimax-robust spectral characteristics of the optimal estimates of the
functionals are proposed for some specific classes of admissible spectral densities.
\end{abstract}

\keywords{Stationary stochastic sequence; optimal linear mean-square estimate; mean square error; least favorable spectral
density; minimax-robust estimate; minimax-robust spectral characteristic}
\maketitle
\noindent{\bf AMS 2010 subject classifications.} Primary: 60G10, 60G25, 60G35, Secondary: 62M20, 93E10, 93E11
\vskip 8pt%
\noindent{\bf DOI:} 10.19139/soic.v4i1.172

\section{Introduction}\label{sc:introduction}
A great number of scientific papers are dedicated to investigation of properties of stationary stochastic processes. Extensive expositions of the theory of stationary stochastic processes can be found, for example, in the books by Gikhman and Skorokhod \cite{Gihman:Skorohod}, Hannan \cite{Hannan}, Rozanov \cite{Rozanov}, Yaglom \cite{Yaglom:1987a,Yaglom:1987b}.
The developed theory of stationary processes has a lot of theoretical and practical applications in engineering, econometrics, finance etc. See, for example, the book by
Box, Jenkins and Reinsel \cite{Box:Jenkins}.
An important problem of the theory is the problem of estimation (extrapolation, interpolation and filtering) of the unobserved values of a stationary process based on observations of the process or on observations of another process which is stationary related with the process under investigation. Constructive methods of extrapolation, interpolation and filtering of stationary processes were developed by Kolmogorov \cite{Kolmogorov}, Wiener \cite{Wiener}, Yaglom \cite{Yaglom:1987a,Yaglom:1987b}. For some recent results see papers by Bondon~\cite{Bondon1},\cite{Bondon2},
Kasahara, Pourahmadi and Inoue~\cite{Kasahara}, \cite{Pourah}, Frank and Klotz~\cite{Frank}, Klotz and M\"adler~\cite{Klotz},
Lindquist \cite{Lindquist}.

The crucial assumption of most of the methods developed for estimation of the unobserved values of stochastic processes is that the spectral densities of the involved stochastic processes are exactly known. However, in practice complete information on the spectral densities is impossible in most cases. In this situation one finds parametric or nonparametric estimate of the unknown spectral density and then apply one of the traditional estimation methods provided that the estimated density is the true one. This procedure can  result in significant increasing of the value of error as
Vastola and Poor~\cite{Vastola} have demonstrated with the help of some examples.
To avoid this effect one can search estimates which are optimal for all densities from a certain class of admissible spectral densities. These estimates are called minimax since they minimize the maximum value of the error.
Grenander \cite{Grenander} was the first one who applied this method to the extrapolation problem for stationary processes.
Franke \cite{Franke:1985}, Franke and Poor \cite{Franke}
investigated the minimax interpolation and extrapolation problem for stationary sequences by using convex optimization techniques.
For more details see a survey paper by Kassam and Poor \cite{Kassam}.

In the papers by Moklyachuk \cite{Moklyachuk:1991} - \cite{Moklyachuk:2008} extrapolation, interpolation and filtering problems for stochastic processes and sequences are investigated.
The corresponding problems for vector-valued stationary sequences and processes are investigated by Moklyachuk and Masyutka~\cite{Moklyachuk:2012}.
In the article by Dubovets'ka and Moklyachuk \cite{Dubovetska4}, \cite{Dubovetska8}  and in the book  by Golichenko and Moklyachuk \cite{Golichenko} the minimax estimation problems are investigated for another generalization of stationary processes --  periodically correlated stochastic sequences and stochastic processes.

Particularly, in the paper
\cite{Moklyachuk:1991} the minimax-robust filtering problem for stationary stochastic sequences from observations with noise is considered. The problem is solved under the condition that the spectral densities $f(\lambda) $ and $f(\lambda)+g(\lambda)$  admit the canonical factorizations. An estimate is found as a solution to the minimization problem by using convex optimization methods. Formulas are proposed that determine the least favorable spectral densities and the minimax-robust spectral characteristics of the optimal estimates
for some special classes of admissible spectral densities.

In this paper we deal with the problem of the mean-square optimal linear estimation of the functionals
$A {\xi}=\sum_{k=0}^{\infty} a(k)\xi(-k)$ and $A_N
{\xi}=\sum_{k=0}^{N} a(k)\xi(-k)$ which depend on the unknown values of a stationary stochastic sequence
$\xi(k)$ with the spectral density $f(\lambda)$ from observations of the sequence $\xi(k)+\eta(k)$ at points of time $k=0,-1,-2,\ldots$, where
$ \eta(k)$ is an uncorrelated with the sequence $\xi(k)$ stationary stochastic sequence with the spectral density $g(\lambda)$.
In contrast to the paper \cite{Moklyachuk:1991}, in this  paper we propose solution of the problem using the method of projection in the Hilbert space of random variables with finite second moments proposed by Kolmogorov \cite{Kolmogorov}. The derived spectral characteristics are determined with the help of operators in the space $\ell_2$ constructed from the Fourier coefficients of the functions
$\frac{1}{f(\lambda)+g(\lambda)}$, $\frac{f(\lambda)}{f(\lambda)+g(\lambda)}$ and $\frac{f(\lambda)g(\lambda)}{f(\lambda)+g(\lambda)}$.  In the case of spectral
uncertainty, where the spectral densities are not exactly
known, but a set of admissible spectral densities is specified, formulas for determination the least favorable spectral densities and the
minimax-robust spectral characteristics of the optimal estimates of the
functionals $A {\xi}$ and $A_N{\xi}$ are proposed for some specific classes of admissible spectral densities.
Such an approach to  filtering problems for stochastic sequence with stationary $n$th increments was applied in the papers by Luz and Moklyachuk \cite{luz3}, \cite{luz6}. The method of projection in the Hilbert space is also applied to extrapolation and interpolation problems for stochastic sequences with stationary $n$th
increments in papers by Luz and Moklyachuk \cite{luz2}, \cite{luz8}.

\section{Classical filtering problem}\label{sc:filtering}

Consider a stationary stochastic sequence $\{\xi(m):m\in\mr Z\}$ with an absolutely
continuous spectral function $F(\lambda)$ which has the spectral density function
$f(\lambda)$.
 Let  $\{\eta(m):m\in\mr Z\}$ be a stationary stochastic sequence which is uncorrelated with $\{\xi(m):m\in\mr Z\}$   and has an absolutely
continuous spectral function $G(\lambda)$ with the spectral density function
$g(\lambda)$. The stationary stochastic sequences $\{\xi(m):m\in\mr Z\}$  and  $\{\eta(m):m\in\mr Z\}$  admit the following spectral representations \cite{Gihman:Skorohod}:
\[\xi (m )=\ip
e^{im\lambda} dZ_{\xi}(\lambda),\quad \eta (m )=\ip
e^{im\lambda} dZ_{\eta}(\lambda),\]
where $Z_{\xi}(\lambda)$ and $Z_{\eta}(\lambda)$ are two random processes with uncorrelated increments defined on $[-\pi,\pi)$ that correspond to the spectral functions $F(\lambda)$ and $G(\lambda)$.

We will suppose that the spectral density functions  $f(\lambda)$ and $g(\lambda)$
satisfy the minimality condition
\be
 \quad\ip
\dfrac{1}{  f(\lambda)+g(\lambda) }
d\lambda<\infty.\label{umova11 fs}\ee
Suppose that we have observations of the sequence
$\xi(m)+\eta(m)$ at points of time  $m=0,-1,-2,\ldots$. The problem is to find the mean-square optimal linear estimates of the functionals
\[A\xi=\sum_{k=0}^{\infty}a(k)\xi(-k)=\ip
A(e^{i\lambda}) dZ_{\xi}(\lambda), \quad A(e^{i\lambda })=\sum_{k=0}^{\infty}a(k)e^{-i\lambda k},\]
\[ A_N\xi=\sum_{k=0}^{N}a(k)\xi(-k)=\ip
A_N(e^{i\lambda}) dZ_{\xi}(\lambda), \quad A_N(e^{i\lambda })=\sum_{k=0}^{N}a(k)e^{-i\lambda k},\]
which depend on the unknown values of the sequence $\xi(k)$, $k\leq0$.

Define by $\widehat{A} \xi$ and $\widehat{A}_N \xi$ the mean-square optimal linear estimates of the values of the functionals $A\xi$ and
$A_N\xi$  based on
observations of the sequence
$\xi(m)+\eta(m)$ at points of time $m=0,-1,-2,\ldots$.
The mean-square errors of the estimates $\widehat{A}\xi$ and $\widehat{A}_N \xi$ are determined as $\Delta(f,g,\widehat{A}\xi):=\mathsf E |A\xi-\widehat{A}\xi|^2$ and $\Delta(f,g,\widehat{A}_N\xi):=\mathsf E |A_N\xi-\widehat{A}_N\xi|^2$.

To find the mean-square optimal estimate of the functional $ A\xi$  we use the method of projection in Hilbert space proposed by Kolmogorov
\cite{Kolmogorov}.

Suppose that the coefficients $\{a(k):k\geq0\}$ defining the functional $A\xi$ satisfy the following conditions:
\be\label{umova na b
fs}\sum_{k=0}^{\infty}|a(k)|<\infty,\quad
\sum_{k=0}^{\infty}(k+1)|a(k)|^2<\infty.\ee
The first condition ensures that the functional $A\xi$ has a finite second moment. The second condition ensures that a linear operator $\me A$ in the space $\ell_2$  to be defined below is compact.

Denote by $H^{0}(\xi +\eta )$
 the closed linear subspace generated by elements
$\{\xi(m)+\eta(m):m=0,-1,-2,\ldots\}$ in the space
$H=L_2(\Omega,\mathcal{F}, \mt P)$ of random variables with $0$ mathematical expectations, $\mathsf E\gamma=0$, and finite variations,  $\mathsf E|\gamma|^2<\infty$.
Let us also consider a subspace $L_2^{0}(f+g)$ in the space $L_2(f+g)$, which is generated by functions
$\{ e^{im\lambda} :m=0,-1,-2,\ldots\}.$

The estimate $\widehat{A} \xi$ of the functional $A\xi$ has the represention
\be \label{otsinka B fs} \widehat{A} \xi=\ip
h (\lambda)dZ_{\xi+\eta }(\lambda), \ee
 where $h (\lambda)$ is the spectral characteristic of the estimate $\widehat{A} \xi$.
The optimal estimate $\widehat{A} \xi$ is a projection of the element $A \xi$ of the space $H$ on the subspace $H^{0}(\xi+\eta)$. It is determined by the following conditions:

1) $ \widehat{A} \xi\in H^{0}(\xi +\eta  ) $;

2) $(A\xi-\widehat{A} \xi)\perp
H^{0}(\xi +\eta )$.

Condition 2) implies that the spectral characteristic $h (\lambda)$ of the optimal estimate $\widehat{A} \xi$ satisfy the following equalities for all $k\leq0$
\[\mathsf E(A\xi-\widehat{A}\xi)
(\overline{\xi (k )+\eta (k )})=\frac{1}{2\pi}\ip
\ld(A(e^{i\lambda
})  -h (\lambda)\rd)
e^{-i\lambda
k} f(\lambda)d\lambda-\frac{1}{2\pi}\ip h (\lambda)e^{-i\lambda
k}  g(\lambda)d\lambda=0
,\]
which leads to the condition
\[\ip\ld(A (e^{i\lambda })
 f(\lambda) -h (\lambda)(f(\lambda)+g(\lambda))\rd)e^{-i\lambda
k} d\lambda=0,\quad \forall k\leq0.\]
It follows from the last equality that the spectral characteristic of the estimate $\widehat{A}\xi$ satisfies the equality
\be \label{spch1}
A (e^{i\lambda })
 f(\lambda) -h (\lambda)(f(\lambda)+g(\lambda))=C (e^{i\lambda}),
 \ee
\[C (e^{i\lambda})=\sum_{k=0}^{\infty}c (k)e^{i\lambda
(k+1)},
\]
where $\{c(k):k\geq0\}$ are coefficients we have to find. Equality (\ref{spch1}) implies that the spectral characteristic of the estimate $\widehat{A}\xi$
can be presented in the following form
\[h
 (\lambda)=A(e^{i\lambda }) \frac{f(\lambda)}{f(\lambda)+g(\lambda)}-\frac{
C (e^{i\lambda})}{ f(\lambda)+g(\lambda) }.\]

Let us find the coefficients $\{c(k):k\geq0\}$.
 Condition 1) implies that the Fourier coefficients of the function $h (\lambda)$ equal to zero if $k>0$. Therefore the following relations hold true:
 \[\ip \ld(A(e^{i\lambda
})\frac{f(\lambda)}{f(\lambda)+g(\lambda)}-
\frac{ C (e^{i\lambda})}{ f(\lambda)+g(\lambda)}\rd)e^{-i\lambda
l}d\lambda=0,\quad l\geq 1. \]
We can write these relations in the form
\be  \sum_{m=0}^{\infty}a(m)\ip \frac{e^{-i\lambda (m+l)}f(\lambda)}{f(\lambda)+g(\lambda)}d\lambda-
\sum_{k=0}^{\infty}c (k)\ip\frac{ e^{i\lambda
(k-l+1)}}{ f(\lambda)+g(\lambda)} d\lambda=0,\quad l\geq 1.\label{spivv_um_1}\ee
  Let us introduce the Fourier coefficients of the functions
\[R_{k,j}=\frac{1}{2\pi}\ip
e^{-i\lambda(j+k)}\frac{f(\lambda)}{f(\lambda)+g(\lambda)}d\lambda;\]
\[P_{k,j} =\frac{1}{2\pi}\ip e^{i\lambda (j-k)}
\dfrac{1 }{ f(\lambda)+g(\lambda)}
d\lambda;\] \[Q_{k,j} =\frac{1}{2\pi}\ip
e^{i\lambda(j-k)}\frac{ f(\lambda)g(\lambda)}
{ f(\lambda)+g(\lambda)}d\lambda.\]

Using the introduced notations we can verify that equality $(\ref{spivv_um_1})$ is equivalent to the following system of equations:
\[\sum_{m=0}^{ \infty}R_{l+1,m}a(m)=
\sum_{k=0}^{\infty}P_{l+1,k+1} c  (k),\quad l\geq0,\]
which admits the matrix representation \[\me R\me a=\me
P \me c ,\] where $\me c =(c (0),c (1),c (2), \ldots)'$, $\me a=(a(0),a(1),a(2),\ldots)'$, $\me
P $, $\me R$ are linear operators in the space $\ell_2$  defined by  matrices with coefficients
$(\me P )_{l,k}=P_{l,k} $, $l,k\geq0$, $(\me
R)_{l, m} =R_{l+1,m}$, $l\geq0$, $m\geq0$. Thus,  the coefficients  $c(k)$, $k\geq0$, can be calculated by the formula
\[c (k)=(\me P ^{-1}\me R\me a)_k,\]  where $(\me
P ^{-1}\me R\me a )_k$ is the
 $k$th element of the vector $\me P ^{-1}\me R\me a$.

The derived expressions allow us to conclude that the spectral characteristic $h (\lambda)$
of the optimal estimate $\widehat{A}\xi$ can be calculated by the formula
\be\label{spectr B
fs}h (\lambda)=A(e^{i\lambda})  \frac{f(\lambda)}{f(\lambda)+g(\lambda)}
-\frac{ \sum_{k=0}^{\infty}
(\me P ^{-1}\me R\me a)_k e^{i\lambda
(k+1)}}{ f(\lambda)+g(\lambda)}.\ee
The value of the mean-square
error of the estimate $\widehat{A}\xi$ of the functional ${A}\xi$ can be calculated by the formula
\[\Delta(f,g;\widehat{A}\xi)=\frac{1}{2\pi}\ip
\frac{\ld|A(e^{i\lambda }) g(\lambda)+ \sum_{k=0}^{\infty} (\me
P ^{-1}\me R \me a)_ke^{i\lambda
(k+1)}\rd|^2}{  (f(\lambda)+g(\lambda))^2}f(\lambda)d\lambda\]
\[
+\frac{1}{2\pi}\ip \frac{\ld|A(e^{i\lambda }) f(\lambda)- \sum_{k=0}^{\infty} (\me
P ^{-1}\me R\me a)_ke^{i\lambda
(k+1)}\rd|^2}{  (f(\lambda)+g(\lambda))^2}
g(\lambda)d\lambda\]
\be\label{poh B fs}=\langle \me R\me
 a,\me P ^{-1}\me R\me a\rangle+\langle\me Q \me a,\me
 a\rangle,\ee where $\me Q $ is a linear operator in the space $\ell_2$ defined by the matrix $\me Q$ with coefficients $(\me Q )_{l,k}=Q_{l,k} $, $l,k\geq0$.

We should note that it might be difficult to find the inverse operator  $ \me P ^{-1}$ if the operator $\me P $ is defined using the Fourier coefficients of the function $\dfrac{1}{ f (\lambda)+g (\lambda)}$. However, the operator $ \me P ^{-1}$ can be easily found if the sequence $\xi(m)+\eta(m)$ admits the canonical moving average representation. In this case the function
$f(\lambda)+g(\lambda)$ admits the canonical factorization \cite{Hannan}
\be  \label{fakt2}
\dfrac{1}{ f(\lambda)+g(\lambda) }=
\ld|\sum_{k=0}^{\infty}\psi (k)e^{-i\lambda k}\rd|^2=\ld|\sum_{k=0}^{\infty}\theta (k)e^{-i\lambda k}\rd|^{-2}.\ee

\begin{lemma}\label{lema_fact_0}
Suppose that canonical factorization ($\ref{fakt2}$) holds true. Let  the linear operators
 $ \Psi$ and $ \Theta$ in the space $\ell_2$ be defined by the matrices with elements $( \Psi )_{k,j}=\psi(k-j)$, $( \Theta )_{k,j}=\theta (k-j)$   for $0\leq j\leq k$, $(
\Psi )_{k,j}=0$, $(
\Theta )_{k,j}=0$ for $0\leq k<j$. Then the linear operator $\me P $  in the space $\ell_2$ admits the factorization $\me
P =\Psi'  \overline{\Psi} $ and the inverse operator $\me
P ^{-1}$ admits the factorization $\me
P ^{-1}= \overline{\Theta} \Theta' $.
\end{lemma}

$\mathrm{\mbox{Proof.}}$ Factorization  ($\ref{fakt2}$) implies
\[\dfrac{1}{f(\lambda)+ g(\lambda)}=\sum_{m=-\infty}^{\infty}p (m)e^{i\lambda m}=\ld|\sum_{k=0}^{\infty}\psi(k)e^{-i\lambda k}\rd|^2\]
\[=\sum_{m=-\infty}^{-1}\sum_{k=-m}^{\infty}\psi (k)\overline{\psi}(k+m)e^{i\lambda m}+\sum_{m=0}^{\infty}\sum_{k=0}^{\infty}\psi(k)\overline{\psi} (k+m)e^{i\lambda m}.\]
Thus,  $p (m)=\sum_{k=0}^{\infty}\psi(k)\overline{\psi} (k+m)$, $m\geq0$, and $p (-m)=\overline {p (m)}$, $m\geq0$. In the case $i\geq j$ we have  \[P_{i,j}=p (i-j)=\sum_{l=i}^{\infty}\psi(l-i)\overline{\psi} (l-j)=(\Psi' \overline{\Psi} )_{i,j},\]
 And in the case $i<j$ we have
 \[P_{i,j}=p (i-j)=\overline{p (j-i)}=\sum_{l=j}^{\infty}\overline{\psi} (l-j)\psi(l-i)
 =(\overline{\Psi}' \Psi)_{j,i}=(\Psi' \overline{\Psi} )_{i,j},\]
which proves the factorization $\me
P =\Psi'  \overline{\Psi} $.

Factorization $(\me
P )^{-1}= \overline{\Theta} \Theta' $ comes from the relation $\Psi\Theta =\Theta \Psi=I$ which has to be proved. If follows from factorization ($\ref{fakt2}$) that
\be\label{spivvidn1}1=\ld(\sum_{k=0}^{\infty}\psi(k)e^{-i\lambda k}\rd)\ld(\sum_{k=0}^{\infty}\theta (k)e^{-i\lambda k}\rd)=\sum_{j=0}^{\infty}\ld(\sum_{k=0}^{j}\psi(k)\theta (j-k)\rd)e^{-i\lambda j},\ee
which leads to the following equalities finishing the proof of the lemma:
\[
    \delta_{i,j}=\sum_{k=0}^{i-j}\psi(k)\theta (i-j-k)=\sum_{p=j}^i\theta (i-p)\psi(p-j)=(\Theta \Psi)_{i,j}=(\Psi\Theta )_{i,j}. \quad \Box\]

Summing up the relations described above, we can formulate the following theorem.

\begin{theorem}\label{thm1 fs}
Let the spectral densities $f(\lambda)$ and $g(\lambda)$ of the uncorrelated stationary stochastic sequences  $\{\xi(m),m\in\mr Z\}$ and
$\{\eta(m),m\in\mr Z\}$
satisfy minimality condition (\ref{umova11 fs}). Let conditions (\ref{umova na b fs}) hold true. The optimal linear estimate
$\widehat{A}\xi$ of the functional $A \xi=\sum_{k=0}^{\infty}a(k)\xi(-k)$ which depends on the unknown values of the stationary sequence $\xi (k )$, $k\leq0$,  based on observations of the sequence $\xi (m )+\eta (m )$
at points of time $m=0,-1,-2,\ldots$ is calculated by the formula
 (\ref{otsinka B fs}).
The spectral characteristic $h (\lambda)$ and the value of the mean-square
error $\Delta(f,g;\widehat{A}\xi)$ of the optimal estimate $\widehat{A}\xi$
 are calculated by formulas  (\ref{spectr B fs}) and (\ref{poh B fs}) respectively.
In the case where the spectral density $f(\lambda)+g(\lambda)$ of the sequence $\xi(m)+\eta(m)$ admits canonical factorization $(\ref{fakt2})$, the operator $\me P^{-1}$ in formulas $(\ref{spectr B
fs})$, $(\ref{poh B fs})$ can be presented as $(\me
P )^{-1}= \overline{\Theta} \Theta'$.
\end{theorem}

To derive the estimate of the functional \[A_N\xi=\sum_{k=0}^{N}a(k)\xi(-k)\] which depends on the unknown values of the stationary stochastic sequence  $\xi(k)$ at points of time $k=0,-1,-2,\ldots,-N$ we can use theorem $\ref{thm1 fs}$. Let us consider the vector $\me a_N=(a(0),a(1),\ldots,a(N), 0,\ldots)'$ to find the spectral characteristic  $h_N (\lambda)$ and the mean-square
error $\Delta(f,g;\widehat{A}_N\xi) $ of the estimate \be \label{otsinka B_N fs} \widehat{A}_N \xi=\ip
h_N (\lambda)dZ_{\xi +\eta }(\lambda) \ee of the functional  $A_N\xi$. In this case the estimate $\widehat{A}_N \xi$ does not depend on the elements $ (\me
R)_{l, m}$ of the operator $\me R$ for $m>N$.
Define a linear operator $\me
R_N$ in the space $\ell_2$   by the matrix with elements $(\me
R_N)_{l, m}=R_{l+1, m}$ for  $l\geq0$, $0\leq m\leq N$, $(\me
R_N)_{l, m}=0$ for $l\geq0$, $m>N$. Then the spectral characteristic of the optimal estimate is calculated by the formula
\be\label{spectr B_N
fs}h_N (\lambda)=A_N(e^{i\lambda})  \frac{f(\lambda)}{f(\lambda)+g(\lambda)}
-\frac{ \sum_{k=0}^{\infty}
(\me P ^{-1}\me R_N\me a_N)_k e^{i\lambda
(k+1)}}{ f(\lambda)+g(\lambda)},\ee
where  $A_N(e^{i\lambda })=\sum_{k=0}^{N}a(k)e^{-i\lambda k}$.
The value of the mean-square error of the optimal estimate $\widehat{A}_N\xi$ is calculated by the formula
\[\Delta(f,g;\widehat{A}_N\xi)=\frac{1}{2\pi}\ip
\frac{\ld|A_N(e^{i\lambda }) g(\lambda)+ \sum_{k=0}^{\infty} (\me
P ^{-1}\me R_N \me a_N)_ke^{i\lambda
(k+1)}\rd|^2}{  (f(\lambda)+g(\lambda))^2}f(\lambda)d\lambda\]
 \[+\frac{1}{2\pi}\ip \frac{\ld|A_N(e^{i\lambda }) f(\lambda)- \sum_{k=0}^{\infty} (\me
P ^{-1}\me R_N\me a_N)_ke^{i\lambda
(k+1)}\rd|^2}{  (f(\lambda)+g(\lambda))^2}
g(\lambda)d\lambda\]
\be\label{poh B_N fs}=\langle \me R_N\me
 a_N,\me P ^{-1}\me R_N\me a_N\rangle+\langle\me Q_N \me a_N,\me
 a_N\rangle,\ee
where $\me Q_N $ is a linear operator  in the space $\ell_2$ which is defined by the matrix with elements $(\me Q_N )_{l,k}=Q_{l,k} $ if  $ 0\leq  l,k\leq N$, and $(\me Q_N )_{l,k}=0 $ if  $ l,k> N$.

The following theorem holds true.

\begin{theorem}\label{thm2 fs}
Let the spectral densities $f(\lambda)$ and $g(\lambda)$ of the uncorrelated stationary stochastic sequences  $\{\xi(m),m\in\mr Z\}$ and
$\{\eta(m),m\in\mr Z\}$ satisfy minimality condition (\ref{umova11 fs}).
The optimal linear estimate
$\widehat{A}_N\xi$  of the functional  $A_N\xi$ which depends on the unknown values $\xi (-k )$, $0\leq k\leq N$,  based on observations of the sequence $\xi (m )+\eta (m )$
at points of time $m=0,-1,-2,\ldots$ is calculated by
formula  (\ref{otsinka B_N fs}).
The spectral characteristic $h_N (\lambda)$ and the value of the mean-square
error $\Delta(f,g;\widehat{A}_N\xi)$ of the estimate $\widehat{A}_N\xi$  are calculated by formulas  (\ref{spectr B_N
fs}) and (\ref{poh B_N fs}).
In the case where the spectral density $f(\lambda)+g(\lambda)$ of the sequence $\xi(m)+\eta(m)$ admits canonical factorization $(\ref{fakt2})$, the operator $\me P^{-1}$ in formulas $(\ref{spectr B_N
fs})$, $(\ref{poh B_N fs})$  can be presented as $(\me
P )^{-1}= \overline{\Theta} \Theta' $.
\end{theorem}

Theorem  $\ref{thm2 fs}$ allows us to find an estimate of the unknown
$\xi (p)$ at point
$p$, $p\leq0$, based on observations of the sequence $\xi(m)+\eta(m)$ at points of time $m=0,-1,-2,\ldots$. Consider the vector $\me a_N$
with $1$ at the coordinate  $(- p) $ and $0$ at the rest of coordinates, and put it in $(\ref{spectr B
fs})$. We get that the spectral characteristic
$\varphi_m(\lambda )$ of the optimal estimate
\be\label{otsinka xi fs} \widehat{\xi} (p)=\ip\varphi_p(\lambda )
dZ_{\xi +\eta }(\lambda)\ee can be presented by the formula
\be\label{spectr xi fs}\varphi_p(\lambda )=e^{i\lambda p}
\frac{f(\lambda)} {f(\lambda)+g(\lambda)}-\frac{ \sum_{k=0}^{\infty}
(\me P ^{-1}\me r_p)_k e^{i\lambda (k+1)}}{  f(\lambda)+g(\lambda
)},\ee where $\me r_p =(R_{0,-p},R_{1,-p},\ldots)$. The value of the mean-square error of the estimate $\widehat{\xi} (p)$ can be calculated by the formula
\[\Delta(f,g;\widehat{\xi} (p))=\frac{1}{2\pi}\ip
\frac{\ld|e^{i\lambda p} g(\lambda)+ \sum_{k=0}^{\infty} (\me
P^{-1}\me r_p )_ke^{i\lambda
(k+1)}\rd|^2}{ (f(\lambda)+g(\lambda))^2}f(\lambda)d\lambda\]
\[
+\frac{1}{2\pi}\ip \frac{\ld|e^{i\lambda
p} f(\lambda)- \sum_{k=0}^{\infty}
(\me P^{-1}\me r_p)_ke^{i\lambda
(k+1)}\rd|^2}{ (f(\lambda)+g(\lambda))^2}
g(\lambda)d\lambda\]
\be \label{pohybka xi fs}=\langle \me r_p,\me P ^{-1}\me r_p\rangle+Q_{-p,-p}.\ee

Thus, the following corollary holds true.

\begin{corollary}
The optimal linear estimate $\widehat{\xi} (p)$ of the unknown value  $\xi (p)$, $p\leq0$,   based on
observations of the sequence $\xi(m)+\eta(m)$ at points of time
$m=0,-1,-2,\ldots$ is defined by formula (\ref{otsinka xi fs}).
The spectral characteristic $\varphi_p(\lambda )$
of the optimal estimate $\widehat{\xi} (p)$ is calculated by  formula
(\ref{spectr xi fs}). The value of the mean-square error
$\Delta(f,g;\widehat{\xi} (p))$ is calculated by  formula
(\ref{pohybka xi fs}).
In the case where the spectral density $f(\lambda)+g(\lambda)$ of the sequence $\xi(m)+\eta(m)$ admits canonical factorization $(\ref{fakt2})$, the operator $\me P^{-1}$ in formulas $(\ref{spectr B
fs})$, $(\ref{poh B fs})$  can be presented as $(\me
P )^{-1}= \overline{\Theta} \Theta' $.
\end{corollary}

 Let us find the mean-square optimal estimate $\widehat{\xi} (0 )$ of the value $ {\xi} (0 )$ of the stochastic sequence  $ {\xi}(m)$  based on
observations of the sequence $\xi(m)+\eta(m)$ at points of time
$m=0,-1,-2,\ldots$. This problem is known as the smoothing of stochastic sequence.

Put $r(k)=R_{k,0}$, $k\in \mr Z$. Then $\{r(k):k\in \mr Z\}$ are the Fourier coefficients of the function $\dfrac{f(\lambda)}{f(\lambda)+g(\lambda)}$. They satisfy the property $r
(k)=\overline{r}(-k)$, $k\in \mr Z$, where  $\overline{r}(k)$ is the complex conjugate to  $r(k)$. Denote by $\{V_{k,j} :k,j\geq0\}$ elements of the matrix which define the operator $\me V=(\me P)^{-1}$. Then
\be\label{oberneni_oper}\sum_{l\geq0}V_{l,j} P_{k,l}=\delta_{kj},\quad k,j\geq0,\ee where $\delta_{kj}$ is the Kronecker delta.  It follows from relations  $(\ref{spectr xi fs})$ and  $(\ref{oberneni_oper})$ that the spectral characteristic $\varphi(\lambda )$ of the estimate $\widehat{\xi} (0 )$ of the element ${\xi} (0 )$
is of the form
 \[\varphi(\lambda )=
\sum_{k=0}^{ \infty}(\overline{r}(k)-(\me Y \me V \me r_0)_k)e^{ -i\lambda k},\]
where $\me Y$ is a linear operator in the space $\ell_2$  defined by the matrix with elements $(\me Y )_{k,l}=P_{k,-l} $, $k\geq0$, $l\geq0$. $(\me Y \me V \me r_0)_k$ is the
 $k$th element of the vector  $\me Y \me V \me r_0$.
The optimal estimate of the element  $ {\xi} (0 )$ is calculated by  the formula \be \label{otsinka xi0 fs} \widehat{\xi} (0 )= \sum_{k=0}^{\infty}(\overline{r}(k)-(\me Y \me V \me r_0)_k)(\xi (-k )+\eta (-k ))
 .\ee
The value of the mean-square error of the estimate is calculated by the formula \be\label{pohybka xi0 fs}\Delta(f,g;\widehat{\xi} (0 ))=\sum_{j=0}^{\infty}\sum_{k=0}^{\infty}
\overline{V}_{k,j} \overline{r}(j)r(k)+\sum_{l\in\mr Z}r(l) g  (-l), \ee
where $\{{g} (k):k\in \mr Z\}$ are the Fourier coefficients of the spectral density $   g(\lambda)   $.

 \begin{corollary}\label{nas xi0 fs}
The optimal linear estimate $\widehat{\xi}(0)$ of the value   $\xi (0 )$ of the stationary stochastic sequence $\xi(m)$ based on observations of the sequence $\xi(m)+\eta(m)$ at points of time
$m=0,-1,-2,\ldots$ is calculated by  formula
(\ref{otsinka xi0 fs}). The value of the mean-square error
$\Delta(f,g;\widehat{\xi} (0 ))$ of the estimate $\widehat{\xi}(0)$ is calculated by  formula
(\ref{pohybka xi0 fs}).
\end{corollary}

Formulas (\ref{spectr B
fs}) and (\ref{poh B fs}) for calculating the spectral characteristic and the value of the mean-square error of the estimate $\widehat{A}\xi$ of the functional  $A\xi$ have been found by the method of orthogonal projection in Hilbert space. Let us show that these formulas are equivalent to the formulas obtained in the paper \cite{Moklyachuk:1991}.
Suppose that canonical factorization (\ref{fakt2}) holds true as well as  the canonical factorization
\be  \label{fakt3}
    f(\lambda)=\sum_{k=-\infty}^{\infty}f(k)e^{i\lambda k}=
    \ld|\sum_{k=0}^{\infty}\phi(k)e^{-i\lambda k}\rd|^2.\ee
Let $\me G$ be a linear operator  in the space $\ell_2$ determined by the matrix with elements $(\me G)_{l,k}=f(l-k)$, $l,k\geq0$.

\begin{lemma}\label{lema_fact_1}
 Suppose the functions $(f(\lambda)+g(\lambda))^{-1}$ and $f(\lambda)$ admit canonical factorizations ($\ref{fakt2}$) and ($\ref{fakt3}$) respectively. Let the linear operators
 $\Psi$ and $\Phi$ in the space $\ell_2$ be determined by  matrices with elements
 $(\Psi)_{k,j}=\psi(k-j)$ and $(
\Phi)_{k,j}=\phi(k-j)$ for $0\leq j\leq k$, $(\Psi)_{k,j}=0$ and $ (\Phi)_{k,j}=0$ for $0\leq k<j$. Then

a) the function $\dfrac{f(\lambda)}{ f(\lambda)+ g(\lambda)}$ admits the canonical factorization
\be  \label{fakt4}
    \dfrac{f(\lambda) }{ f(\lambda)+ g(\lambda)}=\sum_{k=-\infty}^{\infty}s (k)e^{i\lambda k}=
    \ld|\sum_{k=0}^{\infty}\upsilon(k)e^{-i\lambda k}\rd|^2,\ee
where
\[
    \upsilon(k)=\sum_{j=0}^k\psi(j)\phi(k-j)=\sum_{j=0}^k\phi(j)\psi(k-j);\]

b) the linear operator $\Upsilon$ in the space $\ell_2$  defined by the matrix with elements
 $(\Upsilon)_{k,j}=\upsilon(k-j)$ for $0\leq j\leq k$, $(
\Upsilon)_{k,j}=0$ for $0\leq k<j$ can be  presented as
$\Upsilon=\Psi\Phi=\Phi\Psi$.
\end{lemma}

$\mathrm{\mbox{Proof.}}$ Statement a) follows from the equalities
\[
    \ld(\sum_{k=0}^{\infty}\psi(k)e^{-i\lambda k}\rd)\ld(\sum_{k=0}^{\infty}\phi(k)e^{-i\lambda k}\rd)
    =\sum_{j=0}^{\infty}\sum_{k=j}^{\infty}\psi(j)\phi(k-j)e^{-i\lambda k}
    =\sum_{k=0}^{\infty}\ld(\sum_{j=0}^{k}\psi(j)\phi(k-j)\rd)e^{-i\lambda k}.\]
Statement b) follows from the equalities holding true for $i\geq j$:
\[\upsilon(i-j)=\sum_{k=0}^{i-j}\psi(k)\phi(i-j-k)
=\sum_{p=j}^i\phi(i-p)\psi(p-j)=(\Phi\Psi)_{i,j}=(\Psi\Phi)_{i,j}.\quad \Box\]

\begin{lemma}\label{lema_fact_2}
Suppose that canonical factorizations ($\ref{fakt2}$), ($\ref{fakt3}$) hold true and  let the linear operators
  $\Phi$ and $\Upsilon$ be  defined in the same way as in lemma $\ref{lema_fact_1}$. Define  the linear operator $\me T $ in the space $\ell_2$ by the matrix with elements $(\me T )_{l,k}=s (l-k)$, $l,k\geq0$, where the coefficients $s (k)$, $k\geq0$, are defined in (\ref{fakt4}). Then
 operators  $\me T $ and $\me G$ in the space $\ell_2$ admit the factorizations $\me
T =\Upsilon'  \overline{\Upsilon} $ and $\me G=\Phi' \overline{\Phi}$.
\end{lemma}

$\mathrm{\mbox{Proof.}}$ The proof of the lemma is the same as the proof of factorization of the operator  $\me P $ in lemma \ref{lema_fact_0}.
$\Box$

 Lemmas $\ref{lema_fact_0}$, $\ref{lema_fact_1}$ and $\ref{lema_fact_2}$ can be applied to modify  formulas (\ref{spectr B
fs}) and (\ref{poh B fs}).  Put $ \me e =\Theta' \me R \me a $. Canonical factorization ($\ref{fakt2}$) implies
\[
    \frac{
    \sum_{k=0}^{\infty}
    (\me P ^{-1}\me R \me a )_k e^{i\lambda
    (k+1)}}{f(\lambda)+g(\lambda)}=\ld(\sum_{k=0}^{\infty}\psi(k)e^{-i\lambda k}\rd)\sum_{j=0}^{\infty}\sum_{k=0}^{\infty}\overline{\psi} (j)(\overline{\Theta} \me e )_ke^{i\lambda(k+j+1)}\]
\[
    =\ld(\sum_{k=0}^{\infty}\psi(k)e^{-i\lambda k}\rd)\sum_{m=0}^{\infty}\sum_{p=0}^{m}
    \sum_{k=p}^m\overline{\psi} (m-k)\overline{\theta} (k-p)e (p)e^{i\lambda (m+1)}
    =\ld(\sum_{k=0}^{\infty}\psi(k)e^{-i\lambda k}\rd)\sum_{m=0}^{\infty}e (m)e^{i\lambda(m+1)},\]
where $ e (m)=(\Theta' \me R \me a )_m$, $m\geq0$, is the $m$th element of the vector
$\me e =\Theta' \me R \me a $. Since
\[
    (\Theta' \me R \me a )_m=\sum_{j=0}^{\infty}\sum_{p=m}^{\infty}\theta (p-m)s (p+j+1)a(j)=\sum_{j=0}^{\infty}\sum_{l=0}^{\infty}\theta (l)s (m+j+l+1)a(j),\]
the following equality holds true:
\be\label{simple_sp_char_part1_f_st.n_d}
    \frac{
    C (e^{i\lambda})}{ f(\lambda)+g(\lambda)}=\ld(\sum_{k=0}^{\infty}\psi(k)e^{-i\lambda k}\rd)\sum_{m=1}^{\infty}\sum_{j=0}^{\infty}\sum_{l=0}^{\infty}\theta (l)s (m+j+l)a(j)e^{i\lambda m}.\ee

Factorization  ($\ref{fakt4}$) and relation ($\ref{spivvidn1}$) let us make the following transformations:
\[
    \frac{A(e^{i\lambda }) f(\lambda)}{ f(\lambda)+g(\lambda)}
    =\ld(\sum_{k=0}^{\infty}\psi(k)e^{-i\lambda k}\rd)\ld(\sum_{k=0}^{\infty}\theta (k)e^{-i\lambda k}\rd)\sum_{j=0}^{\infty}\sum_{m=-\infty}^{\infty}s (m+j)a(j)e^{i\lambda m}\]
\be\label{simple_sp_char_part2_f_st.n_d}
    =\ld(\sum_{k=0}^{\infty}\psi(k)e^{-i\lambda k}\rd)\sum_{m=-\infty}^{\infty}\sum_{j=0}^{\infty}\sum_{l=0}^{\infty}s (m+j+l)\theta (l)a(j)e^{i\lambda m}.\ee
  Relations ($\ref{simple_sp_char_part1_f_st.n_d}$) and ($\ref{simple_sp_char_part2_f_st.n_d}$) allow us to write the following formula for calculating  the spectral characteristic $h (\lambda)$ of the optimal estimate $\widehat{A}\xi$  using coefficients from canonical factorizations  ($\ref{fakt2}$), ($\ref{fakt3}$):
\[h (\lambda)=\ld(\frac{A(e^{i\lambda }) f(\lambda)}{ f(\lambda)+g(\lambda)}-\frac{
\sum_{k=0}^{\infty}
(\me P ^{-1}\me R \me a )_k e^{i\lambda
(k+1)}}{ f(\lambda)+g(\lambda)}\rd)\]
\[
    =\ld(\sum_{k=0}^{\infty}\psi(k)e^{-i\lambda k}\rd)\sum_{m=0}^{\infty}\sum_{j=0}^{\infty}\sum_{l=0}^{\infty}s (j+l-m)\theta (l)a(j)e^{-i\lambda m}\]
\[
    =\ld(\sum_{k=0}^{\infty}\psi(k)e^{-i\lambda k}\rd)\sum_{m=0}^{\infty}\sum_{j=0}^{\infty}\sum_{l=0}^{\infty}\overline{s} (m-j-l)\theta (l)a(j)e^{-i\lambda m}\]
\be\label{simple_spectr A_f_st.n_d}
    =\ld(\sum_{k=0}^{\infty}\psi(k)e^{-i\lambda k}\rd)\sum_{m=0}^{\infty}( \overline{\me T} \Theta  \me a)_m e^{-i\lambda m}
    =\ld(\sum_{k=0}^{\infty}\psi(k)e^{-i\lambda k}\rd)\sum_{m=0}^{\infty}( \me C\overline{\psi} )_m e^{-i\lambda m}.
\ee
The element  $( \me C\overline{\psi} )_m$, $m\geq0$, from the last relation is the $m$th element of the vector $ \me C \overline{\psi} =\overline{\Psi}' \overline{\me G}\me a$,  $\psi=(\psi(0), \psi(1), \psi(2),\ldots)'$, $\me C$ is a linear operator  defined by the matrix with elements $(\me C)_{k,j}=\me c(k+j)$, $k,j\geq0$, $\me c=\overline{\me G}\me a$, $\me G$ is a linear operator   defined by the matrix with elements $(\me G)_{k,j}= f(k-j)$, $k,j\geq0$. Lemma \ref{lema_fact_2} provides the representation $\me G=\Phi'\overline{\Phi}$ of the operator $\me G$,
where $\Phi$ is a linear operator  defined by the matrix with elements $(\Phi)_{k,j}= \phi(k-j)$, $k,j\geq0$.

The value of the mean-square error $\Delta(f,g;\widehat{A}\xi)$ of the estimate $\widehat{A}\xi$ is calculated by the  formula
\[
    \Delta(f,g;\widehat{A}\xi)
    =\frac{1}{2\pi}\int_{-\pi}^{\pi}|A(e^{i\lambda})|^2f(\lambda)d\lambda+
\frac{1}{2\pi}\int_{-\pi}^{\pi}|h (e^{i\lambda})|^2(f(\lambda)+ g(\lambda))d\lambda\]
\[
    -\frac{1}{2\pi}\int_{-\pi}^{\pi}h (e^{i\lambda})\overline{A(e^{i\lambda})}f(\lambda)d\lambda-
    \frac{1}{2\pi}\int_{-\pi}^{\pi}\overline{h (e^{i\lambda})}A(e^{i\lambda})f(\lambda)d\lambda\]
\be\label{simple_poh A_f_st.n_d}
    =\langle\me G\me a,\me a\rangle
    -\langle\me C  \overline{\psi} ,\me C  \overline{\psi} \rangle.
\ee

In such a way we justify an approach to estimate of the functional  $A\xi$ which uses the method of orthogonal projection in Hilbert space to prove the following result which was obtained in  \cite{Moklyachuk:1991} by direct minimizing the value of the mean-square error of an arbitrary estimate of the functional  $A\xi$.

\begin{theorem}\label{thm3 fs}
Let the stationary stochastic sequences $\{\xi(m),m\in\mr Z\}$ and
$\{\eta(m),m\in\mr Z\}$ have the spectral densities $f(\lambda)$ and $g(\lambda)$ which satisfy minimality condition (\ref{umova11 fs}) and admit canonical factorizations (\ref{fakt2}), (\ref{fakt3}). Suppose also that conditions (\ref{umova na b fs}) hold true. Then the spectral characteristic $h (\lambda)$ of the mean-square optimal estimate $\widehat{A}\xi$
 can be calculated by  formula (\ref{simple_spectr A_f_st.n_d}). The value of the mean-square error $\Delta(f,g;\widehat{A}\xi)$
can be calculated by  formula (\ref{simple_poh A_f_st.n_d}).
\end{theorem}

\begin{example}
Consider two uncorrelated moving average sequences  $\{\xi(m):m\in \mr Z\}$ and $\{\eta(m):m\in \mr Z\}$ with the spectral densities \[f(\lambda)=|1-\phi e^{-i\lambda}|^2\quad\mbox{and}\quad g(\lambda)=|1-\psi e^{-i\lambda}|^2\] respectively, where $|\phi|<1$, $|\psi|<1$.
Consider the problem of finding the mean-square optimal linear estimate
of the functional
$A_1\xi=a\xi(0)+b\xi(-1)$  which depends on the unknown values $\xi(0)$, $\xi(-1)$ based on observations of the sequence $\xi(m)+\eta(m)$ at points of time $m=0,-1,-2,\ldots$. The spectral density of the stochastic sequence  $\{\xi(m)+\eta(m):m\in \mr Z\}$ is of the form
\[f(\lambda)+g(\lambda)=|1-\phi e^{-i\lambda}|^2+|1-\psi e^{-i\lambda}|^2=x|1-ye^{-i\lambda}|^2,\]
 where
\[x=\frac{1}{2}\ld( 2+\phi^2+\psi^2\mp\sqrt{(2+\phi^2+\psi^2)^2-4(\phi+\psi)^2}\rd)\]\[=\frac{1}{2}\ld( 2+\phi^2+\psi^2\mp\sqrt{(\phi^2+\psi^2)^2+4(1-\phi\psi)}\rd);\]
\[y=\frac{1}{2(\phi+\psi)}\ld( 2+\phi^2+\psi^2\pm\sqrt{(\phi^2+\psi^2)^2+4(1-\phi\psi)}\rd).\]
Minimality condition  $(\ref{umova11 fs})$ holds true if $|y|<1$. To determine the spectral characteristic of the optimal estimate
of the functional  $A_1 \xi$ we use formula $(\ref{spectr B_N
fs})$. The operator $\me P$ is defined by the coefficients $(\me P)_{l,k}=\dfrac{y^p}{x(1-y^2)}$, where $p=|k-l|$, $l,k\geq0$. The inverse operator $\me V= \me P^{-1}$ is defined by the coefficients
$(\me V)_{0,0}=x$, $(\me V)_{l,l}=x(1+y^2)$ for $l\geq1$,  $(\me V)_{l,k}=-xy$ for $|l-k|=1$, $l,k\geq0$, and $(\me V)_{l,k}=0$ in other cases. The operator $\me R_1$ is defined by the coefficients
$(\me R_1)_{l,0}=\dfrac{y^{l}(y-\phi)(1-\phi y)}{x(1-y^2)}$ for $l\geq0$, $(\me R_1)_{l,1}=\dfrac{y^{l+1}(y-\phi)(1-\phi y)}{x(1-y^2)}$ for $l\geq0$, and $(\me R)_{l,k}=0$ for $l\geq0$, $k\geq2$.  Let the spectral characteristic of the estimate $\widehat{A}_1\xi$ be presented in the form $h_{1}(\lambda)= \sum _{k=0}^{\infty}w(k)e^{-i\lambda k}
$. Then the optimal estimate of the functional  $A_1\xi$ is
\[
    \widehat{A}_1\xi=\sum_{k=0 }^{\infty }w(k)(\xi(-k)+\eta(-k)),\]
 where
 \[
    w(0)= x^{-1}(a(1-\phi y+\phi^2)+b(y-\phi)(1-\phi y)),\]
 \[
    w(1)=a x^{-1}(y-\phi)(1-\phi y)+b x^{-1}(1-y^2)^{-1}(1-2\phi y+\phi^2-y^3),\]
\[
    w(k)=x^{-1}y^{k-2}(y-\phi)(1-\phi y)(by^2+ay+b)  ,  \quad k\geq2.\]
The value of the mean-square error is
\[
    \Delta(f,g;\widehat{A}_1\xi)=x^{-2}(y-\phi)^2(1-\phi y)^2(1-y^2)^{-1}\]
\[
    +x^{-1}(a^2+b^2)(1+\phi^2\psi^2+(y-\phi-\psi-y\phi\psi)(1-y^2)^{-1})\]\[+2ab x^{-1}(1-y^2)^{-1}((y-\phi-\psi)(1-(\phi+\psi)y+\phi\psi(1+y^2))+y\phi^2\psi^2) .\]

If we take $a=0$, $b=1$ then we will have
\[
    \widehat{\xi}(-1)=\sum_{k=0}^{\infty }w_{-1}(k)(\xi(-k)+\eta(-k)),\]
 where
 \[
    w_{-1}(0)= x^{-1}(y-\phi)(1-\phi y),\quad
    w_{-1}(1)=x^{-1}(1-y^2)^{-1}(1-2\phi y+\phi^2-y^3),\]
\[
    w_{-1}(k)=x^{-1}y^{k-2}(y-\phi)(1-\phi y)(y^2+1),  \quad k\geq2.\]
\[
    \Delta(f,g;\widehat{\xi}(-1))=x^{-2}(y-\phi)^2(1-\phi y)^2(1-y^2)^{-1}
    +x^{-1}(1+\phi^2\psi^2+(y-\phi-\psi-y\phi\psi)(1-y^2)^{-1}).\]

Putting $a=1$, $b=0$ we obtain  solution to the smoothing problem:
\[
    \widehat{\xi}(0)=\sum_{k=0}^{\infty }w_0(k)(\xi(-k)+\eta(-k)),\]
 where
 \[
    w_0(0)= x^{-1}(1-\phi y+\phi^2),\quad
    w_0(k)=x^{-1}y^{k-1}(y-\phi)(1-\phi y),  \, k\geq1.\]
\[
    \Delta(f,g;\widehat{\xi}(0))=x^{-2}(y-\phi)^2(1-\phi y)^2(1-y^2)^{-1}
    +x^{-1}(1+\phi^2\psi^2+(y-\phi-\psi-y\phi\psi)(1-y^2)^{-1}).\]
\end{example}

\section{Minimax-robust method of filtering}\label{sc:Minimax}

  The value of the mean-square error
$\Delta(h (f,g);f,g):=\Delta(f,g;\widehat{A}\xi) $ and
the spectral characteristic $h(\lambda) $
of the optimal linear estimate $\widehat{A}\xi$ of the functional  ${A}\xi$ which depends on the unknown values of the sequence
$\xi(k)$ based on observations of the sequence $\xi(m)+\eta(m)$ are calculated by formulas
$(\ref{spectr
B fs})$ and $(\ref{poh B fs})$ if the spectral densities
$f(\lambda)$ and $g(\lambda)$ of the stochastic sequences $\xi(m)$ and
$\eta(m)$ are exactly known. In the case where a set $\md D=\md
D_f\times\md D_g$ of admissible
spectral densities is given, it is reasonable to apply the minimax-robust method of filtering
of the functionals which consists in  minimizing
the value of the mean-square error for all spectral densities from the
given class $\md D=\md
D_f\times\md D_g$. For description of the proposed method we propose the following definitions \cite{Moklyachuk:1991}.

\begin{definition} For a given class of spectral densities $\mathcal{D}=\md
D_f\times\md D_g$ the spectral densities
$f_0(\lambda)\in\mathcal{D}_f$, $g_0(\lambda)\in\md D_g$ are called
least favorable in the class $\mathcal{D}$ for the optimal linear
filtering of the functional $A\xi$ if the following relation holds
true
\[\Delta(f^0,g^0)=\Delta(h (f^0,g^0);f^0,g^0)=
\max_{(f,g)\in\mathcal{D}_f\times\md
D_g}\Delta(h (f,g);f,g).\]\end{definition}

\begin{definition} For a given class of spectral densities $\mathcal{D}=\md
D_f\times\md D_g$ the spectral characteristic $h^0(\lambda)$ of
the optimal linear estimate of the functional $A \xi$ is called
minimax-robust if there are satisfied conditions
\[h^0(\lambda)\in H_{\mathcal{D}}=\bigcap_{(f,g)\in\mathcal{D}_f\times\md D_g}L_2^{0 }(f+g),\]
\[\min_{h\in H_{\mathcal{D}}}\max_{(f,g)\in
\mathcal{D}_f\times\md D_g}\Delta(h ;f,g)=\sup_{(f,g)\in\mathcal{D}_f\times\md
D_g}\Delta(h^0;f,g).\]\end{definition}

From the introduced definitions and formulas derived in the previous section
we can conclude that the following statement holds true.

\begin{lemma} Spectral densities $f^0\in\mathcal{D}_f$,
$g^0\in\mathcal{D}_g$ satisfying condition $(\ref{umova11
fs})$ are the least favorable in the class $\md D=\md D_f\times\md D_g$ for
the optimal linear filtering of the functional $A\xi$ if  operators $\me
P ^0$, $\me R^0$, $\me Q^0$ determined by the
Fourier coefficients of the functions $\dfrac{ 1}{ f^0 (\lambda)+g^0 (\lambda)}$,
$\dfrac{f^0 (\lambda)}{f^0 (\lambda)+g^0 (\lambda)}$
and $\dfrac{ f^0 (\lambda)g^0 (\lambda)}{  f^0 (\lambda)+g^0 (\lambda) }$
determine a solution to the constrain optimization problem \be \max_{f\in\mathcal{D}}(\langle
\me R  \me a,\me P ^{-1}\me R  \me
a\rangle+\langle\me Q   \me a,  \me a\rangle)= \langle \me R^0  \me a,(\me
P ^0)^{-1}\me R^0  \me a \rangle+\langle\me
Q ^0 \me a,  \me a\rangle. \label{minimax1
fs}\ee
The minimax-robust spectral characteristic is determined as $h^0=h (f^0,g^0)$ if  $h (f^0,g^0)\in H_{\mathcal{D}}$.
\end{lemma}

The minimax-robust spectral characteristic $h^0$ and the pair $(f^0,g^0)$ of the least favorable spectral densities form a saddle point of
the function
$\Delta(h;f,g)$ on the set $H_{\mathcal{D}}\times\mathcal{D}$.
The saddle point inequalities \[\Delta(h;f^0,g^0)\geq\Delta(h^0;f^0,g^0)\geq\Delta(h^0;f,g)
\quad\forall f\in \mathcal{D}_f,\forall g\in \mathcal{D}_g,\forall h\in H_{\mathcal{D}}\]
 hold true if
$h^0=h (f^0,g^0)$ and $h (f^0,g^0)\in H_{\mathcal{D}}$,
where $(f^0,g^0)$ is a solution to the constrain optimization problem
\be \label{extr1}
\widetilde{\Delta}(f,g)=-\Delta(h (f^0,g^0);f,g)\to \inf,\quad (f,g)\in \mathcal{D},
\ee
\[
\Delta(h (f^0,g^0);f,g)
=\frac{1}{2\pi}\ip
\frac{\ld|A (e^{i\lambda }) g^0(\lambda)+ \sum_{k=0}^{\infty} ((\me
P ^0)^{-1}\me R^0 \me a)_ke^{i\lambda
(k+1)}\rd|^2}{ (f^0(\lambda)+g^0(\lambda))^2}f(\lambda)d\lambda
\]
\[+\frac{1}{2\pi}\ip
\frac{\ld|A (e^{i\lambda }) f^0(\lambda)- \sum_{k=0}^{\infty} ((\me
P ^0)^{-1}\me R^0 \me a)_ke^{i\lambda
(k+1)}\rd|^2}{  (f^0(\lambda)+g^0(\lambda))^2}g(\lambda)d\lambda
.\]
This problem is equivalent to the following optimization problem:
\[\Delta_{\mathcal{D}}(f,g)=
\widetilde{\Delta}(f,g)+\delta(f,g|\mathcal{D}_f\times  D_g)\to\inf.\]
A solution of this problem is determined by the condition $0\in
\partial\Delta_{\mathcal{D}}(f^0,g^0)$, which is the necessary and sufficient condition that the pair  $(f^0,g^0)$ belongs to the set of minimums of  the convex functional  $\Delta_{\mathcal{D}}(f,g)$ \cite{Moklyachuk:2008nonsm}, \cite{Pshenychn}, \cite{Rockafellar}. Here   $\partial\Delta_{\mathcal{D}}(f^0,g^0)$ is a subdifferential of the functional  $\Delta_{\mathcal{D}}(f,g)$ at the point $(f,g)=(f^0,g^0)$, which is a set of all linear bounded functionals  $\Lambda$ on $\mathcal{L}_1\times \mathcal{L}_1$ satisfying the inequality
\[\Delta_{\mathcal{D}}(f,g)-\Delta_{\mathcal{D}}(f^0,g^0) \geq\Lambda\ld((f,g)-(f^0,g^0)\rd),\quad  (f,g)\in \mathcal{D}.\]

The form of the functional $\Delta_{\mathcal{D}}(f,g)$ allows us to find derivatives and differentials  in the space $\mathcal{L}_1\times \mathcal{L}_1$ Hence, the complexity of the optimization problem (\ref{extr1}) is determined by the complexity of finding a subdifferential of the indicator function  $\delta(f,g|\mathcal{D}_f\times
    \mathcal{D}_g)$   of the set  $\mathcal{D}_f\times
    \mathcal{D}_g$.

\section{Least favorable densities in the class $\md D_f\times\md D_g$}\label{sc:minimax1}

Consider the problem of minimax-robust estimation of the functional  $A\xi$ based on observations of the sequence $\xi(m)+\eta(m)$ at points of time $m=0,-1,-2,\ldots$ under the condition that the spectral densities $f(\lambda)$ and $g(\lambda)$ belongs to the  set of admissible spectral densities  $\md D=\md
D_f\times\md D_g$, where
  \[\md D_f^0=\ld\{f(\lambda)\mid\frac{1}{2\pi}
\ip
f(\lambda)d\lambda\leq P_1\rd\},\quad\md
D_g^0=\ld\{g(\lambda)\mid\frac{1}{2\pi}\ip g(\lambda)d\lambda\leq
P_2\rd\}.\]

Assume that the spectral densities $f^0 \in\md D_f$, $g ^0\in\md
D_g$ and functions
\be h_{ f}(f ^0,g ^0)=\frac{\ld|A (e^{i\lambda
}) g^0(\lambda)+ \sum_{k=0}^{\infty} ((\me
P ^0)^{-1}\me R^0  \me a)_ke^{i\lambda
(k+1)}\rd|}{  f^0(\lambda)+g^0(\lambda)},\label{hf_filtr_fs}\ee
 \be h_{ g}(f ^0,g ^0)=\frac{\ld|A (e^{i\lambda }) f^0(\lambda)- \sum_{k=0}^{\infty} ((\me
P ^0)^{-1}\me R^0 \me a)_ke^{i\lambda
(k+1)}\rd|}{f^0(\lambda)+g^0(\lambda)}\label{hg_filtr_fs}\ee are bounded. Then the functional
$\Delta(h (f^0 ,g^0 );f,g)$ is continuous and bounded in the space $\md L_1\times\md L_1$.
  To solve the constrain optimization problem (\ref{extr1}) we use  the method of Lagrange multipliers. As a result we obtain the following relations determining the least favorable spectral densities $f^0\in\md D^0_f$, $g^0\in\md D^0_g$:
\be \ld|A (e^{i\lambda }) g^0(\lambda)+ \sum_{k=0}^{\infty} ((\me
P ^0)^{-1}\me R^0  \me a)_ke^{i\lambda
(k+1)}\rd|=\alpha_1  (f^0 (\lambda)+g^0 (\lambda)),\label{D1
rivn1 fs}\ee
\be \ld|A (e^{i\lambda }) f^0(\lambda)- \sum_{k=0}^{\infty} (\me (\me
P ^0)^{-1}\me R^0 \me a)_ke^{i\lambda
(k+1)}\rd|=\alpha_2 (f^0 (\lambda)+g^0 (\lambda)),\label{D1
rivn2 fs}\ee
 where constants $\alpha_1\geq0$ and $\alpha_2\geq0$.
 Particularly, $\alpha_1\neq0$ if $\ip f^0 (\lambda)d\lambda=2\pi P_1$ and $\alpha_2\neq0$ if $\ip g^0 (\lambda)d\lambda=2\pi P_2$. The described relations allow us to formulate the following theorems.

\begin{theorem} Suppose that the spectral densities $f^0(\lambda)\in\md D_f$ and
$g^0(\lambda) \in\md D_g$ satisfy minimality condition $(\ref{umova11
 fs})$ and the functions $h_{ f}(f ^0,g ^0)$ and
$h_{ g}(f ^0,g ^0)$ determined by formulas (\ref{hf_filtr_fs}) and (\ref{hg_filtr_fs}) are bounded. The functions
$f^0 (\lambda)$ and $g^0 (\lambda)$ determined by equations
$(\ref{D1 rivn1 fs})$, $(\ref{D1 rivn2 fs})$ are the least favorable in the class $\md D=\md D_f \times\md D_g$ for the optimal linear filtering of the functional  $A\xi$
 if they determine a solution to
optimization problem $(\ref{minimax1 fs})$.
 The function  $h (f ^0,g ^0)$ determined by formula  $(\ref{spectr B
fs})$
 is the minimax-robust spectral characteristic of the optimal estimate of the functional  $A\xi$.\end{theorem}

\begin{theorem} Suppose that the spectral density $f(\lambda) $ is known, the spectral density $g^0(\lambda)\in\md D_g$ and they  satisfy minimality condition
$(\ref{umova11 fs})$. Suppose also that the function
$h_{ g}(f ,g ^0)$ determined  by formula (\ref{hg_filtr_fs}) is bounded. The spectral density  $g ^0(\lambda)$ is the least favorable in the class $\md D_g$ for the optimal linear filtering of the functional $A\xi$ if it is of the form
\[g ^0(\lambda)=\max\ld\{0,\frac{\ld|A (e^{i\lambda }) f (\lambda)- \sum_{k=0}^{\infty} ((\me
P ^0)^{-1}\me R^0  \me a)_ke^{i\lambda
(k+1)}\rd|}{\alpha_2 }-f(\lambda)\rd\}\]
 and the pair $(f ,g ^0)$ is a solution of the
optimization problem $(\ref{minimax1 fs})$. The function $h (f ,g ^0)$ determined by formula
$(\ref{spectr B fs})$  is the minimax-robust spectral characteristic of the optimal estimate of the
functional  $A\xi$.\end{theorem}

\section{ Least favorable densities in the class $\md D_{f,g}$}\label{sc:minimax2}

Consider the problem of optimal linear filtering of the functional
$A\xi$ in the case, where the  spectral densities
$f(\lambda)$ and $g(\lambda)$ belong to the set of admissible spectral densities $\md D_{f,g}$:
 \[\md D_{f,g}=\ld\{(f(\lambda), g(\lambda))\mid\frac{1}{2\pi}
\ip
\frac{1}{f(\lambda)+g(\lambda)}d\lambda\geq P_0\rd\}.\]

Suppose that the pair of spectral densities $(f^0,g^0)\in \md D_{f,g}$ is such that the functions
$h_{ f}(f ^0,g ^0)$,
$h_{ g}(f ^0,g ^0)$ determined by $(\ref{hf_filtr_fs})$, $(\ref{hg_filtr_fs})$ are bounded.
     As in the previous section, we use  the method of Lagrange multipliers to solve the constrain optimization problem (\ref{extr1}). As a result we obtain the following relations determining the least favorable spectral densities  $f^0\in\md D^0_f$, $g^0\in\md D^0_g$:
 \be\ld|A (e^{i\lambda }) g^0(\lambda)+ \sum_{k=0}^{\infty} ( (\me
P ^0)^{-1}\me R^0  \me a)_ke^{i\lambda
(k+1)}\rd|^2=\beta_1^2,\label{D3 rivn1
fs}\ee
  \be\ld|A (e^{i\lambda }) f^0(\lambda)- \sum_{k=0}^{\infty} (\me (\me
P ^0)^{-1}\me R^0 \me a)_ke^{i\lambda
(k+1)}\rd|^2=\beta_2^2,\label{D3 rivn2 fs}\ee where
$\beta_1\geq0$ and $\beta_2\geq0$.
Suppose that the least favorable spectral densities  $f^0(\lambda)$ and $g^0(\lambda)$ admits a representation \[f^0(\lambda)=\sum_{k=-\infty}^{\infty}f^0(k)e^{i\lambda k},\quad g^0(\lambda)=\sum_{k=-\infty}^{\infty}g^0(k)e^{i\lambda k}.\]
Denote by $\me A$ the operator determined by the matrix with elements $(\me A)_{k,j}=a(k+j)$, $k,j\geq0$, and denote by $\me A^+$ the operator determined by the matrix with elements $(\me A^+)_{k,j}=a(j-k)$ for $j\geq k\geq0$, $(\me A^+)_{k,j}=0$ otherwise. Denote also by $\me f^0$ and  $\me g^0$  vectors  determined by coefficients $(\me f^0)_0=\dfrac{f^0(0)}{2}$, $(\me f^0)_j=f^0(j)$, $j\geq1$, and $(\me g^0)_0=\dfrac{g^0(0)}{2}$, $(\me g^0)_j=g^0(j)$, $j\geq1$. Then equations $(\ref{D3 rivn1
fs})$ and $(\ref{D3 rivn2
fs})$ can be written as
\begin{equation}\label{D4 rivn1
fs}\ld|\sum_{k=0}^{\infty}\ld(\ld(\me A^+\rd)'\me g^0+\me A\me g^0\rd)_ke^{-i\lambda k}+\sum_{k=0}^{\infty} \ld(\me A^+\me g^0+ (\me
P ^0)^{-1}\me R^0  \me a\rd)_ke^{i\lambda
(k+1)}\rd|^2=\beta_1^2
\end{equation}
\begin{equation}\label{D4 rivn2
fs}\ld|\sum_{k=0}^{\infty}\ld(\ld(\me A^+\rd)'\me f^0+\me A\me f^0\rd)_ke^{-i\lambda k}-\sum_{k=0}^{\infty}\ld(\me A^+\me f^0+ (\me
P ^0)^{-1}\me R^0  \me a\rd)_ke^{i\lambda
(k+1)}\rd|^2=\beta_2^2
\end{equation}

\begin{theorem} Suppose that a pair of spectral densities $(f ^0, g^0)\in\md D_{f,g}$ satisfy minimality  condition $(\ref{umova11
 fs})$ and the functions $h_{ f}(f ^0,g ^0)$,
$h_{ g}(f ^0,g ^0)$ determined by formulas (\ref{hf_filtr_fs}) and (\ref{hg_filtr_fs}) are bounded. The functions
$f^0 (\lambda)$ and $g^0 (\lambda)$ determined by equations
$(\ref{D4 rivn1 fs})$ and $(\ref{D4 rivn2 fs})$ are the least favorable in the class $\md D_{f,g}$ for the optimal linear filtering of the functional  $A\xi$
if they determine a solution to
optimization  problem $(\ref{minimax1 fs})$.
The function  $h (f ^0,g ^0)$ determined by formula  $(\ref{spectr B
fs})$
 is the minimax-robust spectral characteristic of the optimal estimate of the
functional  $A\xi$.
\end{theorem}

\section{ Least favorable densities in the class $\md D=\md D^u_v\times\md D_{\eps}$}\label{sc:minimax3}

Consider the problem of optimal linear filtering of the functional
$A\xi$ for the  set of admissible spectral densities  $\md D=\md
D^u_v\times\md D_{\eps}$, \[\md
D^u_v=\ld\{f(\lambda)\mid v(\lambda)\leq f(\lambda)\leq
u(\lambda),\,\frac{1}{2\pi}\ip f(\lambda)d\lambda\leq P_1\rd\},\]
\[\md
D_{\eps}=\ld\{g(\lambda)\mid g(\lambda)=(1-\eps)g_1(\lambda)+\eps
w(\lambda),\frac{1}{2\pi}\ip g(\lambda)d\lambda\leq P_2\rd\},\]
where
the spectral densities $u(\lambda)$, $v(\lambda)$, $g_1(\lambda)$
are known and fixed  and the spectral densities $u(\lambda)$ and $v(\lambda)$
are bounded.

Suppose that the spectral densities  $f ^0\in \md D^u_v$, $g ^0\in \md D_{\eps}$ and the functions $h_{ f}(f ^0,g ^0)$,
$h_{ g}(f ^0,g ^0)$ defined in
$(\ref{hf_filtr_fs})$, $(\ref{hg_filtr_fs})$ are bounded. Then condition $0\in\partial\Delta_{\md
D}(f^0 ,g^0 )$ implies the following equations determining the least favorable spectral densities:
 \be\ld|A (e^{i\lambda }) g^0(\lambda)+ \sum_{k=0}^{\infty} ( (\me
P ^0)^{-1}\me R^0  \me a)_ke^{i\lambda
(k+1)}\rd|=    (f^0 (\lambda)+g^0 (\lambda))
(\gamma_1(\lambda)+\gamma_2(\lambda)+\alpha_1^{-1}),\label{D2 rivn1
fs}\ee
  \be\ld|A (e^{i\lambda }) f^0(\lambda)- \sum_{k=0}^{\infty} (\me (\me
P ^0)^{-1}\me R^0 \me a)_ke^{i\lambda
(k+1)}\rd|= (f^0 (\lambda)+g^0 (\lambda))
(\varphi(\lambda)+\alpha_2^{-1}),\label{D2 rivn2 fs}\ee where
$\gamma_1\leq0$ and $\gamma_1=0$ if $f ^0(\lambda)\geq
v(\lambda)$; $\gamma_2(\lambda)\geq0$ and $\gamma_2=0$ if
$f ^0(\lambda)\leq u(\lambda)$; $\varphi(\lambda)\leq0$ and
$\varphi(\lambda)=0$ if $g ^0(\lambda)\geq (1-\eps)g_1(\lambda)$.
Thus, the following theorem holds true.

\begin{theorem} Suppose that $f^0(\lambda)\in \md D^u_v$, $g^0(\lambda)\in \md
D_{\eps}$ and condition $(\ref{umova11 fs})$ holds true. Let
the functions $h_{ f}(f ^0,g ^0)$ and
$h_{ g}(f ^0,g ^0)$ determined by formulas (\ref{hf_filtr_fs}) and (\ref{hg_filtr_fs}) are bounded.
The functions $f^0 (\lambda)$ and $g^0 (\lambda)$ determined by equations $(\ref{D2 rivn1 fs})$ and $(\ref{D2 rivn2 fs})$ are the least
favorable in the class $\md D=\md D^u_v \times\md D_{\eps}$ for the optimal linear filtering of the functional  $A\xi$
 if they determine
a solution to optimization problem
  $(\ref{minimax1 fs})$.
 The function  $h (f ^0,g ^0)$ determined by formula  $(\ref{spectr B
fs})$
 is the minimax-robust spectral characteristic of the optimal estimate of the
functional
 $A\xi$.\end{theorem}

\begin{theorem} Let the spectral density $f (\lambda)$ be known,
the spectral density $g^0 \in\md D_{\eps}$ and they satisfy minimality condition $(\ref{umova11 fs})$. Assume  that the
function
$h_{ g}(f ,g ^0)$ determined by formula  (\ref{hg_filtr_fs}) is bounded. The spectral density $g ^0(\lambda)$ is the least favorable in the class
  $\md D_{\eps}$ for the optimal linear filtering of the functional $A\xi$ if it is of the form
\[g ^0(\lambda)=\max\ld\{(1-\eps)g_1(\lambda),f_1(\lambda)\rd\},\]
\[f_1(\lambda)= \alpha_2\ld|A (e^{i\lambda }) f (\lambda)- \sum_{k=0}^{\infty} ((\me
P ^0)^{-1}\me R^0 \me a)_ke^{i\lambda
(k+1)}\rd| -f (\lambda)\] and
the pair $(f ,g ^0)$ determines
a solution to optimization problem
$(\ref{minimax1 fs})$. The function  $h (f,g ^0)$ determined by formula  $(\ref{spectr B
fs})$
 is the minimax-robust spectral characteristic of the optimal estimate of the
functional
 $A\xi$.\end{theorem}

\section{Conclusions}\label{sc:Conclusions}

In the article we presented methods of the mean-square optimal linear filtering of the linear functionals $A {\xi}=\sum_{k=0}^{\infty} a(k)\xi(-k)$ and $A_N
{\xi}=\sum_{k=0}^{N} a(k)\xi(-k)$  which depend on the unknown values of a stationary stochastic sequence
$\xi(k)$ with the spectral density $f(\lambda)$ from observations of the sequence $\xi(m)+\eta(m)$ at points of time $m=0,-1,-2,\ldots$, where
$\eta(k)$ is an uncorrelated with the sequence $\xi(k)$ stationary stochastic sequence with the spectral density $g(\lambda)$.
 Formulas for calculating the spectral characteristics and the values of the mean-square errors of the functionals are proposed in the case of spectral certainty, where the spectral densities $f(\lambda)$ and $g(\lambda)$ of the stationary stochastic sequences $\xi(m)$ and $\eta(m)$ are exactly known.
In the case of spectral uncertainty, where the spectral densities are not known, but a set $\md D=\md
D_f\times\md D_g$ of admissible
spectral densities is given, the minimax-robust method of filtering is applied. For some given classes of admissible spectral densities we found relations that determine the least favorable spectral densities and the minimax-robuct spectral characteristics of the estimate of the functionals.


\begin{thebibliography}{99}


\bibitem{Bondon1}
{\sc P. Bondon}, \ \emph{Influence of missing values on the prediction of a stationary time series}, Journal of Time Series Analysis  \textbf{26} (2005), No. 4,  519--525.

\bibitem{Bondon2}
{\sc P. Bondon}, \ \emph{Prediction with incomplete past of a stationary process}. Stochastic Process and Their Applications \textbf{98} (2002), 67--76.


\bibitem{Box:Jenkins}
\newblock G. E. P. Box, G. M. Jenkins and G. C. Reinsel,
\newblock \emph{Time series analysis. Forecasting and control. 3rd ed.},
\newblock Englewood Cliffs, NJ: Prentice Hall, 1994.



\bibitem{Dubovetska4}
\newblock I. I. Dubovets'ka and M. P. Moklyachuk,
\newblock \emph{Filtration of linear functionals of periodically correlated sequences},
\newblock Theory of Probability and Mathematical Statistics, vol. 86, pp. 51-64, 2013.

\bibitem{Dubovetska8}
\newblock I. I. Dubovets'ka and M. P. Moklyachuk,
\newblock \emph{On minimax estimation problems for periodically correlated stochastic processes},
\newblock Contemporary Mathematics and Statistics, vol.2, no. 1, pp. 123-150, 2014.

\bibitem{Frank}
\newblock M. Frank and L. Klotz,
 \newblock \emph{A duality method in prediction theory of multivariate stationary sequences},
 \newblock Math. Nachr. vol.244, pp. 64-77, 2002.


\bibitem{Franke:1985}
\newblock J. Franke,
\newblock \emph{Minimax robust prediction of discrete time series},
\newblock Z. Wahrsch. Verw. Gebiete, vol. 68, pp. 337-364, 1985.

\bibitem{Franke}
\newblock J. Franke and H. V. Poor,
\newblock \emph{Minimax-robust filtering and finite-length robust predictors},
\newblock Robust and Nonlinear Time Series Analysis. Lecture Notes in Statistics, Springer-Verlag,
vol. 26, pp. 87-126, 1984.

\bibitem{Gihman:Skorohod}
\newblock I. I. Gikhman and A. V. Skorokhod,
\newblock \emph{The theory of stochastic processes. I.},
\newblock Berlin: Springer, 2004.

\bibitem{Golichenko}
\newblock I. I. Golichenko and M. P. Moklyachuk,
\newblock \emph{Estimates of functionals of periodically correlated processes},
\newblock Kyiv: NVP ``Interservis", 2014.

\bibitem{Grenander}
\newblock U. Grenander,
\newblock \emph{A prediction problem in game theory},
\newblock Arkiv f\"or Matematik, vol. 3, pp. 371-379, 1957.

\bibitem{Hannan}
\newblock E. J. Hannan,
\newblock \emph{Stationary stochastic processes. 2nd rev. ed.},
\newblock John Wiley \& Sons, New York, 2009.

\bibitem{Karhunen}
\newblock K. Karhunen,
\newblock \emph{Uber lineare Methoden in der Wahrscheinlichkeitsrechnung},
\newblock Annales Academiae Scientiarum Fennicae. Ser. A I, no. 37, 1947.

\bibitem{Kasahara}
\newblock  Y. Kasahara, M. Pourahmadi and A. Inoue,
\newblock \emph{Duals of random vectors and processes with applications to prediction problems with missing values},
\newblock Stat. Probab. Lett. Vol. 79, No. 14, 1637-1646, 2009.

\bibitem{Kassam}
\newblock S.A. Kassam and H. V. Poor,
\newblock \emph{Robust techniques for signal processing: A survey},
\newblock Proceedings of the IEEE, vol. 73, no. 3, pp. 433-481, 1985.

\bibitem{Klotz}
\newblock L. Klotz and C. M\" adler,
 \newblock \emph{On the notion of minimality of a q-variate stationary sequence},
 \newblock Complex Anal. Oper. Theory, 2015

\bibitem{Kolmogorov}
\newblock A. N. Kolmogorov,
\newblock \emph{Selected works by A. N. Kolmogorov. Vol. II: Probability theory and mathematical statistics. Ed. by A.
N. Shiryayev. Mathematics and Its Applications. Soviet Series. 26. Dordrecht etc.}
\newblock Kluwer Academic Publishers, 1992.

\bibitem{Lindquist}
\newblock A. Lindquist and G. Picci,
\newblock \emph{Linear stochastic systems. A geometric approach to modeling, estimation and identification},
\newblock Series in Contemporary Mathematics 1. Berlin: Springer, 2015.

\bibitem{luz2}
\newblock M. M. Luz and M. P. Moklyachuk,
\newblock \emph{Interpolation of functionals of stochastic sequenses with stationary increments},
\newblock Theory of Probability and Mathematical Statistics,  vol. 87, pp. 117-133, 2013.

\bibitem{luz3}
\newblock M. M. Luz and M. P. Moklyachuk,
\newblock \emph{Minimax-robust filtering
problem for stochastic sequence with stationary increments},
\newblock Theory of Probability and Mathematical Statistics,  vol. 89, pp. 127 - 142, 2014.

 \bibitem{luz6}
\newblock M. Luz and M. Moklyachuk,
\newblock \emph{Minimax-robust filtering problem for stochastic sequences with stationary increments and cointegrated sequences},
\newblock Statistics, Optimization \& Information Computing, vol. 2, no. 3, pp. 176 - 199, 2014.


\bibitem{luz8}
\newblock M. Luz and M. Moklyachuk,
\newblock \emph{Minimax-robust prediction problem for stochastic sequences with stationary increments and cointegrated sequences},
\newblock Statistics, Optimization \& Information Computing, vol. 3, no. 2, pp. 160 - 188, 2015.

\bibitem{Moklyachuk:1991}
\newblock M. P. Moklyachuk,
\newblock \emph{Minimax filtering of linear transformations of stationary sequences},
\newblock  Ukrainian Mathematical Journal, vol. 43, no. 1 pp. 92-99, 1991.


\bibitem{Moklyachuk:2000}
\newblock M. P. Moklyachuk,
\newblock \emph{Robust procedures in time series analysis},
\newblock Theory of Stochastic Processes, vol. 6, no. 3-4, pp. 127-147, 2000.

\bibitem{Moklyachuk:2001}
\newblock M. P. Moklyachuk,
\newblock \emph{Game theory and convex optimization methods in robust estimation problems},
\newblock Theory of Stochastic Processes, vol. 7, no. 1-2, pp. 253-264, 2001.

\bibitem{Moklyachuk:2008}
\newblock M. P. Moklyachuk,
\newblock \emph{Robust estimations of  functionals of stochastic processes},
\newblock Kyiv University, Kyiv, 2008.

\bibitem{Moklyachuk:2008nonsm}
\newblock M. P. Moklyachuk,
\newblock \emph{Nonsmooth analysis and optimization},
\newblock Kyiv University, Kyiv, 2008.


\bibitem{Moklyachuk:2012}
\newblock M. Moklyachuk, and O. Masyutka,
\newblock \emph{Minimax-robust estimation technique for stationary stochastic processes},
\newblock LAP LAMBERT Academic Publishing, 2012.

\bibitem{Pourah}
\newblock M. Pourahmadi, A. Inoue and Y. Kasahara
\newblock \emph{A prediction problem in $L^2(w)$}.
\newblock Proceedings of the American Mathematical Society. Vol. 135, No. 4, pp. 1233-1239, 2007.


\bibitem{Pshenychn}
\newblock B. N. Pshenichnyi,
\newblock \emph{Necessary conditions of an extremum},
\newblock ``Nauka'', Moskva, 1982.

\bibitem{Rockafellar}
\newblock R. T. Rockafellar,
\newblock \emph{Convex Analysis},
\newblock Princeton University Press, 1997.


\bibitem{Rozanov}
\newblock Yu. A. Rozanov,
\newblock \emph{Stationary stochastic processes},
\newblock San Francisco-Cambridge-London-Amsterdam: Holden-Day, 1967.


\bibitem{Vastola}
\newblock  K. S. Vastola and H. V. Poor,
\newblock  \emph{ An analysis of the effects of spectral uncertainty on Wiener filtering},
\newblock  Automatica, vol. 28, pp. 289--293, 1983.


\bibitem{Wiener}
\newblock  N. Wiener,
\newblock \emph{Extrapolation, Interpolation and Smoothing of Stationary Time Series. With Engineering Applications},
\newblock  The M. I. T. Press, Massachusetts Institute of Technology, Cambridge, Mass., 1966.

\bibitem{Yaglom:1987a}
\newblock A. M. Yaglom,
\newblock \emph{Correlation theory of stationary and
related random functions. Vol. 1: Basic results},
\newblock Springer Series in Statistics, Springer-Verlag, New York etc., 1987.

\bibitem{Yaglom:1987b}
\newblock A. M. Yaglom,
\newblock \emph{Correlation theory of stationary and
related random functions. Vol. 2: Suplementary notes and references},
\newblock Springer Series in Statistics, Springer-Verlag, New York etc., 1987.



\end{thebibliography}
\end{document}